\documentclass[12pt,a4paper]{amsart}
\usepackage{amssymb,amsmath,amsthm}
\usepackage{graphicx}
\usepackage{xcolor}

\usepackage{epsfig}

\long\def\onefigure#1#2{%  #1 picture,  #2  caption
\begin{figure*}[tbp]
\begin{center}
#1
\end{center}
\caption{#2}
\end{figure*}
}%end onefigure def

\newcommand{\lipefig}[2]  % labeled Ipe figure
{\onefigure{\mbox{\psfig{file=#1.eps}}}{\label{f:#1} #2} }

\newtheorem{theorem}{Theorem}[section]
\newtheorem{lemma}{Lemma}[section]

\newtheorem{claim}{Claim}[section]

\newcommand{\la}{\lambda}

\newcommand{\R}{\mathbb{R}}

\newcommand{\W}{\mathcal{W}}
\newcommand{\X}{\mathcal{X}}
\newcommand{\zz}{\mathcal{Z}}

\newcommand{\inn}{\mathrm{int\,}}
\newcommand{\inc}{\mathrm{int\,conv\,}}
\newcommand{\conv}{\mathrm{conv\;}}

\newcommand{\aff}{\mathrm{aff\,}}

\newcommand{\pos}{\mathrm{pos\,}}
\newcommand{\lin}{\mathrm{lin\,}}

\newcommand{\remove}[1]{}

\numberwithin{equation}{section}

\begin{document}

\title{A point in the interior of the convex hulls}
\author{Imre B\'ar\'any and Yun Qi}
\keywords{Convex hull, theorems of Carath\'eodory and Steinitz}
\subjclass[2020]{Primary 52A20, secondary 52A37}

\maketitle

\begin{abstract}Steinitz's theorem states that if a point $a \in \inc X$ for a set $X \subset \R^d$, then $X$ contains a subset $Y$ of size at most $2d$ such that 
$a \in \inc Y$. The bound $2d$ is best possible here. We prove the colourful version of this theorem and characterize the cases when exactly $2d$ sets are needed. 
\end{abstract}

\section{Introduction and some background}\label{sec:introd}

\bigskip
In a famous paper Steinitz \cite{Ste} proved that if a point $a$ lies in the interior of the convex hull of a set $X \subset \R^d$, i. e. $a \in \inc X$, then $X$ contains a subset $Y$ of size at most $2d$ such that $a \in \inc Y$. The bound $2d$ is best possible as shown by the example $X=\{\pm e_1,\ldots,\pm e_d\}$, where $e_1,\ldots,e_d$ is a basis of $\R^d$. We may assume that $a=0$. Then $a=0\in \inc X$ is the same as $\pos X=\R^d$ where $\pos X$ stands for the positive (or cone) hull of $X$, that is, the set of all linear combinations of elements in $X$ with nonnegative coefficients. 

\smallskip
With $a=0$ we may assume that $0 \notin X$. The condition $\pos X=\R^d$ remains valid if we replace some $x \in X$ by $\la x$ with any $\la >0$. So we may suppose that $X \subset S^{d-1}$, the Euclidean unit sphere of $\R^d$. With this notation Steinitz theorem has the following form.

\begin{theorem}\label{th:stein} If $X \subset S^{d-1}$ and $\pos X=\R^d$, then there is $Y \subset X$ with $\pos Y=\R^d$ and $|Y|\le 2d$. Further there is such a $Y$ with $|Y|\le 2d-1$ unless $X=\{\pm e_1,\ldots,\pm e_d\},$ where $e_1,\ldots,e_d$ is a basis of $\R^d$. 
\end{theorem}

This paper is about a colourful version of Steinitz's result in which there is a system of sets $\X=\{X_1,\ldots,X_{2d}\}$ where each $X_i \subset \R^d$. A {\sl transversal} of $\X$ is a set $T=\{x_1,\ldots,x_{2d}\}$ with $x_i \in X_i$ for every $i\in [2d]:=\{1,\ldots,2d\}$. A {\sl partial transversal} associated with a proper subset $I\subset [2d]$ is a set $T=\{x_i\in X_i: i \in I\}$ and a $k$-transversal is a partial transversal with $|I|=k$.

\begin{theorem}\label{th:colstein} If $X_i \subset S^{d-1}$ and $\pos X_i=\R^d$ for all $i \in [2d]$, then there is a transversal $T$ of the system $\{X_i: i \in [2d]\}$ with $\pos T=\R^d$. 
\end{theorem}

This result is stated informally and without proof in \cite{Bar82} and is attributed to S. Dancs. It also appears in~\cite{Bar21} as Exercise 10.3. Its simple proof is given in Section~\ref{sec:ThmB}. Steinitz's theorem follows from its colourful version by taking $X_i=X$ for all $i \in [2d]$. 

\smallskip
When do we need exactly $2d$ sets in Theorem~\ref{th:colstein}? The main target in this paper is to give a complete answer to this question. 

\smallskip
One case is simple: when $X=\{\pm e_1,\ldots,\pm e_d\}$ and $e_1,\ldots,e_d$ is a basis of $\R^d$ with $\|e_i\|=1$, then the system $X=X_1=\ldots=X_{2d}$ satisfies the condition $\pos X_i=\R^d$ (for every $i$), and it is clear that this system has no $(2d-1)$-transversal $T$ with $\pos T=\R^d$.  So this is one case when exactly $2d$ sets are needed. We call this the {\sl Basis Case} or {\sl BCase} for short. 

\smallskip
This is just a simple variant of the equality case in Theorem~\ref{th:stein}. Surprisingly this is not the only case when equality occurs in the colourful case. There is another one, namely the following.

\begin{figure}[h]
\centering
\includegraphics[scale=1.0]{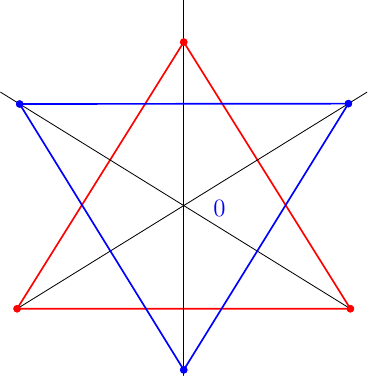}
\caption{The Positive Basis case for $d=2$.}
\label{fig:2D}
\end{figure}

Let $F=\{f_1,\ldots,f_{d+1}\} \subset S^{d-1}$ with $\pos F=\R^d$, $F$ is just the set of vertices of a simplex containing the origin in its interior. Define $X_1=\ldots=X_d=F$ and $X_{d+1}=\ldots=X_{2d}=-F$. Figure~\ref{fig:2D} shows the 2-dimensional case where $X_1=X_2$ is the set of vertices of the blue triangle and $X_3=X_4$ is that of the red one. It is not hard to check (we omit the simple proof) that this system has no  $(2d-1)$-transversal $T$ with $\pos T=\R^d$.  So this is another case when exactly $2d$ sets are needed. We call this one the {\sl Positive Basis Case} or {\sl PCase} for short. 

\smallskip
We note that in BCase there are exactly $(2d)!$ transversals with $\pos T=\R^d$, while in the PCase the number of such transversals is $(d+1)!d!$ as one can see easily. 

\smallskip
Here comes the main result of this paper.

\begin{theorem}\label{th:main}  $2d$ sets are needed in Theorem~\ref{th:colstein} exactly in the BCase or in the PCase.
\end{theorem}

\bigskip
\section{Proof of Theorem~\ref{th:stein}}\label{sec:ThmA}

There are several proofs, see for instance \cite{Ste}, \cite{DGK}, and \cite{Eck}. We present the one below because its method will be used later. The following lemma, the cone version of Carath\'eodory's theorem (see \cite{Cara} and \cite{DGK} or \cite{Bar21}) will be needed.

\begin{lemma}\label{l:Cara} Assume the vector $v\ne 0$ lies in $\pos A$ for some set $A \subset \R^d$. Then there is $B \subset A$ with $v \in \pos B$ and $|B|\le d$.
\end{lemma}

The proof of Theorem~\ref{th:stein} begins by choosing $d$ linearly independent vectors $x_1,\ldots,x_d$ from $X$. Setting $\overline{X}=\{x_1,\ldots,x_d\}$ the cone $C:=\pos \overline{X}$ has nonempty interior and we choose a vector $v \in \inn C$. Then there is a small ball $B(v,\delta)\subset \inn C$ centred at $v$ with radius $\delta>0$. Since $-v \in \pos X$, Lemma \ref{l:Cara} implies that there is a set $\overline{Y}=\{y_1,\ldots,y_k\}\subset X$ with $-v \in \pos \overline{Y}$ and $k\le d$. The origin lies in the interior of $\conv(\{-v\}\cup B(v,\delta))$. Setting $Y=\overline{X}\cup \overline{Y}$ we have $|Y|\le d+k\le 2d$ and $\pos Y=\R^d$. This completes the proof of the first part of Theorem~\ref{th:stein}. 

\smallskip
For the second part assume that $\pos X=\R^d$ and that there is no $Y\subset X$ with $|Y|<2d$ and $\pos Y=\R^d$. Consider $\overline{X}=\{x_1,\ldots,x_d\}$, $v$, and $\overline{Y}=\{y_1,\ldots,y_k\}$ from the previous proof. We see that $k=d$ must hold and $x_i\ne y_j$ for every $i,j\in [d]$ as otherwise $\pos Y=\R^d$ still holds but $|Y|<2d$. Moreover $y_1,\ldots,y_d$ are linearly independent because otherwise $-v$ is in the cone hull of a proper subset $Y'$ of the $\overline{Y}$ and $|\overline{X} \cup Y'|<2d$ and $\pos(\overline{X} \cup Y')=\R^d.$

\begin{figure}[h]
\centering
\includegraphics[scale=0.8]{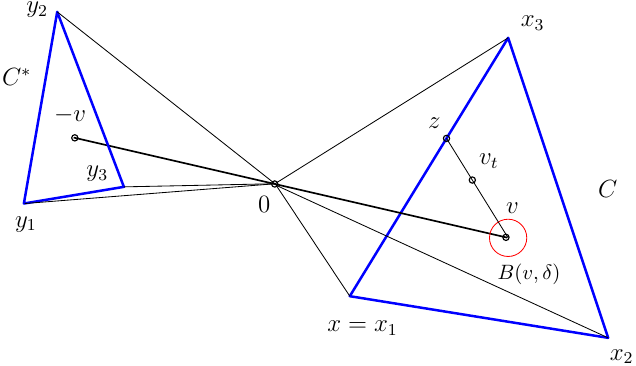}
\caption{ $v \in \inn C$ and $-v \in \inn C^*$.}
\label{fig:Cone}
\end{figure}
\smallskip
Fix now the vector $v$ and the set $Y$. The cone $C^*=\pos \overline{Y}$ has nonempty interior and $-v \in \inn C^*$. Choose a vector $z$ from the boundary of cone $C$, see Figure~\ref{fig:Cone}. Then $v_t=(1-t)v+tz \in \inn C$ for every $t\in [0,1)$ and $-v_t \in \inn C^*$ for small enough $t>0$. We repeat the previous argument with $v_t,-v_t$ (instead of $v,-v$) while $C,C^*$ remain unchanged. We increase $t$ and stop when $v_t \in \inn C$ but $-v_t$ reaches the boundary of $C^*$. In this case a small ball $B(v_t,\delta_t) \subset \inn C$ but $-v_t$ lies in the cone hull of a proper subset $Y'$ of $\overline{Y}$. Thus, the origin lies in the interior of $\conv(\{-v_t\}\cup B(v_t,\delta_t))$ and then $\pos (\overline{X} \cup Y')=\R^d$. A contradiction so $-v_t \in \inn C^*$ for every $t\in [0,1)$.

\smallskip
This implies that $-u \in \inn C^*$ for every $u \in \inn C$. In other words $-\inn C \subset \inn C^* $. The same argument starting with $y_1,\ldots,y_d$ in place of $x_1,\ldots,x_d$ shows that $C=-C^*$, and $X$ contains $\{\pm e_1,\ldots,\pm e_d\}$ for a basis $e_1,\ldots,e_d$ of $\R^d$. Finally one needs to show that $X$ coincides with $\{\pm e_1,\ldots,\pm e_d\}$. We leave the simple proof to our imaginary reader.\qed

\bigskip
\section{The matrix representation of $\X$}\label{sec:matrix}

Before moving on we specify how a transversal of the system $\X$ is given in the BCase and in the PCase. This will only be needed in the last sections.

\smallskip
In the BCase we associate with this system a $2d \times 2d$ matrix $M(\X)=(x_{i,j})$ where $x_{i,j}=e_j$ if $j\le d$ and $x_{i,j}=-e_{j-d}$ if $j>d$. The rows of $M(\X)$ represent the elements in $X_i$. It is evident that a transversal $T$ is then given by a permutation $p(\cdot)$ of $[2d]$, that is, $T$ contains the element $x_{p(h),h}$ from column $h$ of $M(\X)$.

\smallskip
It is clear that $\pos T=\R^d$ for such a transversal $T$. In fact every $b\in \R^d$ can be written uniquely as $b=\sum_1^d \beta_he_h=\sum_{\beta_h>0} \beta_he_h +\sum_{\beta_h<0}\beta_he_h$. If $\beta_h>0$, then $x_{p(h),h}=e_h$ and when $\beta_h<0$, then $x_{p(h+d),h+d}=-e_h$. Consequently $b$ can be written as a non-negative combination of the elements of $T$:
\begin{equation}\label{eq:bpos}
b=\sum_{h \in [d], \beta_h>0} \beta_h x_{p(h),h}+\sum_{h \in [d], \beta_h<0} (-\beta_h)x_{p(h+d),h+d}.
\end{equation}

\smallskip
In the PCase we have a $2d \times(d+1)$ matrix $M(\X)=(x_{i,j})$ where $x_{i,j}=f_j$ if $i\le d$ and $x_{i,j}=-f_j$ if $i>d$. Again, the rows of $M$ represent the elements in $X_i$. We observe that if $T$ is a transversal of $\X$ with $\pos T=\R^d$ and no partial transversal $T'$ of $T$ satisfies $\pos T'=\R^d$, then according to the uniqueness part of Theorem~\ref{th:stein}, $T=\{\pm e_1,\ldots,\pm e_d\}$ where $e_1,\ldots,e_d$ is a basis of $\R^d$. This is only possible if $T=(F\cup -F)\setminus \{f_k,-f_k\}$ for some $k \in [d+1]$.

\smallskip
Then such a transversal $T$ is given by some $k\in [d+1]$ and by a one-to-one map $p^+: ([d+1]\setminus k)\to [d]$ and by another one-to-one map $p^-: ([d+1]\setminus k)\to  \{d+1,\ldots,2d\}$. In this case $T$ contains two elements from column $h$ of $M(\X)$ (for every $h\ne k$), namely $x_{p^+(h),h}$ and $x_{p^-(h),h}$. 

\smallskip 
When $T$ is given by the missing $k\in [d+1]$ and $p^+,p^-$ and $b\in \R^d$, we have a unique linear combination
$b=\sum_{h\ne k} \beta_hf_h=\sum_{\beta_h>0,h\ne k} \beta_hf_h +\sum_{\beta_h<0, h\ne k}(-\beta_h)(-f_h)$
As $f_h=x_{p^+(h),h}$ and $-f_h=x_{p^-(h),h}$ with uniquely determined $p^+(h)$ and $p^-(h)$ this is of the form 
\begin{equation}\label{eq:bPos}
b=\sum_{h: \beta_h>0,h\ne k}\beta_h x_{p^+(h),h}+\sum_{h: \beta_h<0,h\ne k}(-\beta_h)x_{p^-(h),h}.
\end{equation}

\bigskip
\section{Proof of Theorem~\ref{th:colstein}}\label{sec:ThmB}

For this proof we need the colourful version of Lemma~\ref{l:Cara}, see \cite{Bar82} or \cite{Bar21}.

\begin{lemma}\label{l:colCara} Assume the vector $v\ne 0$ lies in $\pos A_i$ for all $A_i \subset \R^d$ $i\in [d]$. Then there is a transversal $T$ of the system $A_i$ with $v \in \pos T$.
\end{lemma}

The proof of Theorem~\ref{th:colstein} is similar to that of Theorem~\ref{th:stein}. 
Let $x_i \in X_i$ for $i\in [d]$ be linearly independent vectors, and define $C=\pos\{x_1,\ldots,x_d\}$. Choose again $v \in \inn C$ so a small ball $B(v,\delta) \subset \inn C$. Since $-v \in \pos X_i$ for $i=d+1,\ldots,2d$ Lemma~\ref{l:colCara} implies the existence of a (possibly partial) transversal $Y'$ of the system $X_{d+1},\ldots,X_{2d}$ with $-v \in \pos Y'$. Of course $|Y'|\le d$. The origin lies in the interior of $\conv(\{-v\}\cup B(v,\delta))$ which is contained in $\pos Y$ where $Y=\{x_1,\ldots,x_d\}\cup Y'$. Thus $\pos Y=\R^d$ and $|Y|\le 2d.$\qed

\bigskip
\section{Preparations for the proof of Theorem~\ref{th:main}}\label{sec:prep}

A set $A$ is a {\sl positive basis} of $\R^d$ if $\pos A =\R^d$ but $\pos (A\setminus \{a\}) \ne \R^d$ for any $a \in A$. One form of Steinitz's theorem is that a positive basis of $\R^d$ has at most $2d$ elements. We will need the following result which is in Davis~\cite{Dav} and in a slightly different form in Reay~\cite{Reay}.

\begin{lemma}\label{l:reay} If $\pos A=\R^d$, then there is a subset $B$ of $A$ that is a positive basis of $\lin B$, the linear span of $B$, and $|B|=\dim \lin B +1$.
\end{lemma}

For the proof of Theorem~\ref{th:main} we use induction on $d$, and the case $d=1$ is trivial. So we assume that $d>1$ and that the statement holds in all dimensions less than $d$. We have a system of sets $X_i \subset \R^d$ with $\pos X_i=\R^d$ for every $i \in [2d]$ and suppose that no $(2d-1)$-transversal $T$ of the system satisfies $\pos T=\R^d$. 
\smallskip
For $v \in X:=\bigcup_1^{2d}X_i$ define $P(v)=\{j \in [2d]: -v \in X_j\}.$ We begin with a simple yet important lemma.

\begin{lemma}\label{l:xand-x} For every $i \in [2d]$ and for every $x \in X_i$ there is a set $J\subset [2d]$ such that $-x \in X_j$ for every $j \in J$, $i\notin J$, and $|J|\ge d$. In other words, $|P(x)\setminus \{i\}|\ge d$ for every $i \in [2d]$ and for every $x \in X_i$.
\end{lemma}

The {\bf proof} goes by an argument similar to the one for Theorem~\ref{th:stein}. For simpler notation we suppose $i=1$. Let $J \subset [2d]$ be the set of all subscripts $j\ne 1$ with $-x \in X_j$. 

\smallskip
\begin{figure}[h]
\centering
\includegraphics[scale=0.8]{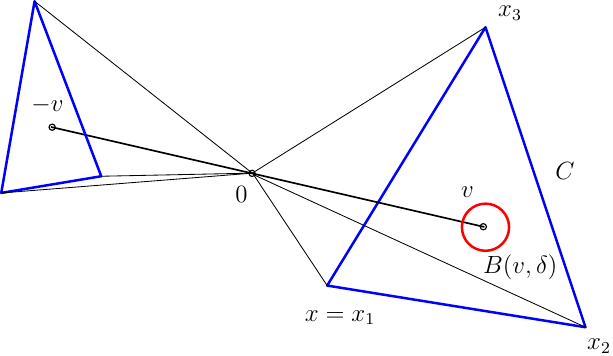}
\caption{$B(v,\delta)\subset\inn C$.}
\label{fig:Ccone}
\end{figure}

We {\bf claim} that $|J|\ge d$. Assume that, on the contrary, $|J| <d$. Choose a set $I\subset [2d]$ containing $J\cup \{1\}$ with $|I|=d$, and vectors $x_i \in X_i$ for all $i \in I$ with $x_1=x \in X_1$ so that these $d$ vectors are linearly independent. This is clearly possible. Set $\overline{X}=\{x_i:i\in I\}$. The interior of the cone $C=\pos \overline{X}$ is nonempty and we pick $v \in \inn C$. Again there is a small ball $B(v,\delta)$ (centred at $v$ and with radius $\delta>0$) contained in $\inn C$. The vector $-v \in \pos X_h$ for every $h \notin I$, see Figure~\ref{fig:Ccone}.

\smallskip
By Lemma~\ref{l:colCara}, there is $x_h \in X_h$ for every $h \notin I$ such that $-v \in \pos \{x_h: h \notin I\}$. Define again $\overline{Y}=\{x_h:h\notin I\}$ and $C^*=\pos \overline{Y}.$ As we have seen in the proof of the first part of Theorem~\ref{th:stein}, $0 \in \inc(\overline{X}\cup \overline{Y})$ and then $\pos Y=\R^d$ where $Y=\overline{X}\cup \overline{Y}$. Moreover, by our assumption, no proper subset of $Y=\{x_1,\ldots,x_{2d}\}$ can have the origin in the interior of its convex hull. As we have seen in the proof of the second part of Theorem~\ref{th:stein}, we have to have $C=-C^*$ and $\overline{X}=-\overline{Y}$. In other words $\{x_i:i \in I\}=\{-x_h:h \notin I\}$. In particular, there is an $h \notin I$ such that $x=x_1=-x_h$ implying that $-x\in X_h$ for some $h \notin J$ (and $h\ne 1$). So $h$ is not in $J$ contradicting the maximality of $J$.\qed

\medskip
In the proof of Theorem~\ref{th:main} we will distinguish two cases: {\bf Case I} when $X_i\cap (-X_i)=\emptyset$ for every $i \in [2d]$, and  {\bf Case II} when $X_i\cap (-X_i)\ne \emptyset$ for some $i \in [2d]$.

\bigskip
\section{Proof of Theorem~\ref{th:main} in Case I}\label{sec:caseI}

In this case $P(v)\cap P(-v) =\emptyset$ for every $v \in X$ and Lemma~\ref{l:xand-x} implies that $P(v)$ and $P(-v)$ form a partition of $[2d]$ for every $v \in X$. As $\pos X_{2d}=\R^d$, $X_{2d}$ contains a positive basis of $\R^d$. This positive basis contains, by Lemma~\ref{l:reay}, a set $\{a_1,\ldots,a_k\}\subset X_{2d}$ which is a positive basis for $L:=\lin \{a_1,\ldots,a_k\}$ whose dimension is $k-1$. Of course $k\le d+1$. Moreover, $k>2$, as $k=2$ would imply $a_1=-a_2$ contrary to $X_{2d}\cap (-X_{2d})=\emptyset$. 

\smallskip
Assume first that $k \le d$. Then by Hall's classic theorem  on distinct representatives~\cite{Hall} (or see for instance \cite{LiW}), the sets $P(a_1),\ldots,P(a_k) \subset [2d]$ have a transversal $i_1 \in P(a_1),\ldots,i_k \in P(a_k)$ with all $i_j$ distinct. This follows from the Hall condition: for every $h\le k$, the union of any $h$ of the sets $P(a_i)$ has size at least $d\ge k\ge h$, as every $|P(a_i)|\ge d$.

\smallskip
For simpler notation we assume that $i_1=1,\ldots,i_k=k$ and then $j=i_j\in P(a_j)$ means that $-a_j \in X_j$ for all $j \in [k]$. The transversal $-a_1,\ldots,-a_k$ for the system $X_1\ldots,X_k$ is a positive basis of $L$. Define the map $x \to  x^*$ as the orthogonal projection of $x \in \R^d$ to the subspace $L^{\perp}$, the orthogonal complementary subspace of $L$. Then $\pos X_j^*=L^{\perp}$ for every $j \in [2d]$ because $\pos X_j=\R^d$. We use Theorem~\ref{th:colstein} in the space $L^{\perp}$ (whose dimension is $d-(k-1)<d$) to the sets $X_{k+1}^*,\ldots, X_{k+2(d-k+1)}^*$; their number is $2(d-k+1)$. So this system has a transversal $x_{k+1}^*\in X_{k+1}^*,\ldots, x_{k+2(d-k+1)}^* \in X_{k+2(d-k+1)}^*$ with $\pos \{x_{k+1}^*\ldots, x_{k+2(d-k+1)}^*\}=L^{\perp}$.

\smallskip
Set $x_i=-a_i$ for $i \in [k]$ and for $i>k$ choose $x_i \in X_i$ whose projection to $L^{\perp}$ is $x_i^*$. This $x_i$ may not be unique in which case any such $x_i$ will do. It is easy to check (we omit the details) that $\pos \{x_1,\ldots, x_{k+2(d-k+1)}\}=\R^d$. This is a transversal $T$ of the system $X_1,\ldots,X_{k+2(d-k+1)}$ with $|T|=k+2(d-k+1)<2d$ because $k>2$, contradicting the original assumption that no such transversal exists. This finishes the case $k\le d$.

\smallskip
Hall's theorem can be used even if $k=d+1$ and gives a system of distinct representatives $i_1,\ldots,i_{d+1}$ for $P(a_1),\ldots,P(a_{d+1})$ as long as $|\bigcup_1^{d+1}P(a_j)| \ge d+1$, all other unions of the sets have size at least $d$, as every $P(a_i)$ is of size at least $d$. Assume this is the case, and take again $i_j=j$ for every $j \in [d+1]$. Here $j\in P(a_j)$ means that $-a_j \in X_j$ for every $j\in [d+1]$ and then $\pos \{a_1,\ldots,a_{d+1}\}=\R^d$. We have a $(d+1)$-transversal $T=\{-a_1,\ldots,-a_{d+1}\}$ of the system $X_1,\ldots,X_{2d}$ with $\pos T=\R^d$, a contradiction again. 

\smallskip
So we are left with the case $|\bigcup_1^{d+1} P(a_j)| = d$ implying that 
\[
P(a_1)=\ldots=P(a_{d+1}) \mbox{ and } P(-a_1)=\ldots=P(-a_{d+1}), 
\]
exactly the Positive Basis Case.\qed

\bigskip
\section{Proof of Theorem~\ref{th:main} in Case II}\label{sec:caseII}

This time there is a vector $v$ and an $i\in [2d]$ with both $v$ and $-v$ in $X_i$. We write $\ell$ for the line $\{\la v: \la \in \R\}$ and $H$ for the $(d-1)$-dimensional subspace $\ell^{\perp}$. Without loss of generality, we may assume that $i=2d$. According to Lemma~\ref{l:xand-x}, $|P(v)\setminus \{2d\}|\ge d$ and $|P(-v)\setminus \{2d\}|\ge d$. So there is another $j\in [2d],\; j\ne 2d$ with $j \in P(v)\cap P(-v)$, that is, both $v$ and $-v$ are in $X_j$. We may assume that $j=2d-1$.

\smallskip
Define $\pi:\R^d \to H$ as the orthogonal projection to $H$ and set $X_i^{\circ}=X_i\setminus \{v,-v\}$ and $Z_i:=\pi(X_i^{\circ})$ for every $i\in [2d]$. Then $\pos Z_i=H$ for every $i\in [2d-2]$, and Theorem~\ref{th:colstein} guarantees the existence of a transversal $Z=\{z_1,\ldots,z_{2d-2}\}$ of the system $\zz=\{Z_1,\ldots,Z_{2d-2}\}$ with $\pos Z=H$. Let $x_i \in X_i$ be a point with $\pi(x_i)=z_i$ chosen arbitrarily if there is more than one such point. We will come back to the uniqueness of $x_i$ in Lemma~\ref{l:uniq}. Then $T=\{x_1,\ldots,x_{2d-2}\}$ is a transversal of the system $\X^{\circ}=\{X_1^{\circ},\ldots,X_{2d-2}^{\circ}\}$. We observe that $\conv T$ is $d$- or $(d-1)$-dimensional because its projection to $H$ is $(d-1)$-dimensional.

\smallskip
If it is $d$-dimensional, then it contains a point $\la v \in \ell$ (for some $\la \in \R$) together with a small (Euclidean) ball $B(\la v,\delta)$ where $\delta>0$. If $\la \ge 0$, then extend $T$ to the $(2d-1)$-transversal (of the system $\X$) $T\cup \{-v\}$ by adding $-v \in X_{2d-1}$ to it. In this case $\pos (T\cup \{-v\})=\R^d$, contrary to our indirect assumption. The same argument works when $\la <0$ with the $(2d-1)$-transversal $T\cup \{v\}$.

\smallskip
We may now assume that $\conv T$ is $(d-1)$-dimensional. In this case $\aff T$ is a hyperplane in $\R^d$ and $ \ell \cap \conv T$ is a single point $\la v \in \ell$ with some $\la \in \R$. Since $\pos Z=H$, $\conv T$ contains $B(\la v,\delta) \cap \aff T$ for some $\delta>0$. If $\la >0$, then with the previous transversal $\pos (T\cup \{-v\})=\R^d$, a contradiction again. The same argument works when $\la <0$. So $\la=0$, $\aff T=\lin T$ is a $(d-1)$-dimensional subspace of $\R^d$ and $\pos T=\lin T$.

\smallskip
Now we return to the choice of $x_i\in X_i$ with $\pi(x_i)=z_i$.

\begin{lemma}\label{l:uniq} The solution $x_i$ to $\pi(x_i)=z_i$ with $x_i\in X_i$ is unique.
\end{lemma}

The {\bf proof} is short. Assume there are $x_i,y_i \in X_i$ with $\pi(x_i)=\pi(y_i)=z_i$. Then $x_i=y_i+\mu v$ with some $\mu \in \R$. We want to show that $\mu=0$. The previous proof goes through without any change when we use $y_i$ instead of $x_i$ and gives the transversal $T(y_i)$ instead of $T=T(x_i)$. Again $\aff T(y_i)=\lin T(y_i)$ is a $(d-1)$-dimensional subspace and $\lin T\cap \lin T(y_i)$ is a $(d-2)$-dimensional subspace unless $\mu = 0$. Moreover, every $x_j$ ($j\ne i$)  lies in this subspace. But then $\pos T$ is only halfspace in $\lin T$. But we just proved that $\pos T= \lin T$. So $\mu=0$.\qed 

\smallskip
Assume next that a partial transversal $Z'$ (which is just $Z$ with one element deleted) also satisfies $\pos Z'=H$ and let $T'$ be the corresponding $(2d-3)$-transversal of $\X^{\circ}$. The previous argument with $T'$ in place of $T$ works again and shows that $\aff T'=\lin T'$ and that $\conv T'$ contains $(\lin T')\cap \rho B$ with some $\rho>0$, where $B$ is the Euclidean unit ball of $\R^d$. The $(2d-1)$-transversal $T^*$ is the extension of $T'$ by adding $v\in X_{2d-1}$ and $-v\in X_{2d}$ to it. Then $\pos T^*=\R^d$, a contradiction again. 

\smallskip
This means that $\pos Z_i=H$ for every $i\in [2d-2]$ and no partial transversal $Z'$ satisfies $\pos Z'=H$. Now we can use induction except that 
the condition $\|z\|=1$ for every $z \in \bigcup_1^{2d-2}Z_i$ does not hold. Yet induction still works for the system $W_i=\{z/\|z\|:z\in Z_i\}$ ($i\in [2d-2]$) as their cone hull is $H$ and they have no partial transversal whose positive (or cone) hull is $H$. By the induction hypothesis the system $\W=\{W_1,\ldots,W_{2d-2}\}$ is either in the BCase or in the PCase.

\smallskip
In the BCase there is a basis $G=\{g_1,\ldots,g_{d-1}\}$ of $H$ such that $W_i=G\cup (-G)$ for every $i\in [2d-2]$. So $|Z_i|=2(d-1)$ and 
then $|X_i^{\circ}|=2(d-1)$ again for every $i$. 

\begin{lemma}\label{l:subBC} There are linearly independent unit vectors $e_1,\ldots,e_{d-1} \in \R^d$ such that $X_i^{\circ}=\{\pm e_1,\ldots,\pm e_{d-1}\}$ for every $i\in [2d-2]$.
\end{lemma}

The proof is not difficult but technical and is postponed to Section~\ref{sec:subBC}

\smallskip
This lemma shows that every $X_i^{\circ}$ lies in a $(d-1)$-subspace $H^{\circ}$, say. But $X_i$ must have points on both sides of $H^{\circ}$ (because $\pos X_i=\R^d$) and the projection of these points must appear in $Z_i$ unless these points coincide with $v$ and $-v$. Thus $v,-v \in X_i$ for all $i \in [2d-2]$, and setting $e_d:=v$ we have $X_i=\{\pm e_1,\ldots,\pm e_d\}$ for all $i \in [2d-2]$. At the start of this proof we had two special sets $X_{2d}$ and $X_{2d-1}$ containing $v$ and $-v$. We see now that we could have chosen any two sets for this special property. It follows that $X_i=\{\pm e_1,\ldots,\pm e_d\}$ for all $i \in [2d]$. This finishes the proof of the BCase.

\smallskip
In the PCase half of the sets in $\W$ coincide with $G=\{g_1,\ldots,g_d\}$ and the other half with $-G$. Here the vectors $g_1,\ldots,g_d$ form a positive basis of $H$. For simpler notation assume $W_i=G$, for $i \in [d-1]$ and $W_i=-G$, for $i \in \{d,\ldots,2d-2\}$. It follows from Lemma~\ref{l:uniq} that $|X_i^{\circ}|=d$ for every $i\in [2d-2]$.

\begin{lemma}\label{l:subPC} There are linearly independent unit vectors $f_1,\ldots,f_{d} \in \R^d$ such that $X_i^{\circ}=F=\{f_1,\ldots,f_d\}$ for $i \in [d-1]$, and $X_i^{\circ}=-F$ for $i \in \{d,\ldots,2d-2\}$ and the linear and cone hull of $F$ is a $(d-1)$-dimensional subspace $H^{\circ}$.
\end{lemma}

The proof is postponed to Section~\ref{sec:subPC}.

\smallskip
We see again that every $X_i^{\circ}$ lies in the $(d-1)$-subspace $H^{\circ}$. The previous argument shows that $v,-v \in X_i$ for every $i\in[2d-2]$. 
Thus $|X_i|=d+2$ for every $i\in [2d-2]$.
As $v,-v \in X_i$ for every $i\in[2d]$ we can repeat this proof with any two sets $X_i,X_j$ instead of $X_{2d}$ and $X_{2d-1}$.
Starting with $X_1$ and $X_2$ shows that $W_{2d-1}=W_{2d}=G$, while starting with $X_{d-1},X_{d}$ yields that 
$W_{2d-1}=-W_{2d}$. A contradiction showing that the PCase cannot occur in case II. \qed

\bigskip
\section{Proof of Lemma~\ref{l:subBC}}\label{sec:subBC}

\bigskip
We are going to use the matrix representations from Section~\ref{sec:matrix}. In the BCase $M(\W)$ is a $(2d-2)\times (2d-2)$ matrix $(w_{i,j})$ with $w_{i,j}=f_j$ if $j<d$ and $w_{i,j}=-f_{j-d+1}$ if $j\ge d$. A transversal $W$ of $\W$ is given by a permutation $p(\cdot)$ of $[2d-2]$ as $W=\{w_{p(1),1},\ldots,w_{p(2d-2),2d-2}\}$ with $w_{p(h),h} \in W_{p(h)}$. Recall that every $x_{i,j} \in X_i^{\circ}$ is mapped to $z_{i,j}=\pi(x_{i,j}) \in Z_i$ which is mapped further to $w_{i,j}=\frac {z_{i,j}}{\|z_{i,j}\|}\in W_i$. Thus 
$z_{i,j}=\|z_{i,j}\|w_{i,j}$. By Lemma~\ref{l:uniq}, $z_{i,j}$ determines $x_{i,j}$ uniquely and $x_{i,j}=z_{i,j}+\lambda_{i,j}v$ where $\lambda_{i,j}\in \R$.

 The transversal $T=\{x_{p(1),1},\ldots,x_{p(2d-2),2d-2}\}$ of $\X^{\circ}$ is defined by the permutation $p$ of $[2d-2]$. Here $x_{p(h),h}\in X_{p(h)}$. We extend $T$ to a transversal $T^*$ of the original system $\X$ by adding $v\in X_{2d-1}$ and $-v \in X_{2d}$ to it.   %the transversal $Z=(z_{p(1),1},\ldots,w_{p(2d-2),2d-2})$$ of $\zz$.

\begin{claim}\label{cl:Pos} $\pos T^*=\R^d$.
\end{claim}

\smallskip
{\bf Proof.} A given $a \in \R^d$ can be written (uniquely) as $a=b+\lambda v$ where $\lambda \in \R.$ Then $b=\sum_1^{d-1}\beta_he_h$, a unique linear combination again. According to equation (\ref{eq:bpos}) (where all sums are taken over $h\in [d-1]$)
\begin{eqnarray*}
b&=&\sum_{\beta_h>0} \beta_h w_{p(h),h}+\sum_{\beta_h<0} (-\beta_h)w_{p(h+d-1),h+d-1}\\
 &=&\sum_{\beta_h>0} \frac{\beta_h}{\|z_{p(h),h}\|} z_{p(h),h}+\sum_{\beta_h<0} \frac{-\beta_h}{\|z_{p(h+d-1),h+d-1}\|}z_{p(h+d-1),h+d-1}\\
 &=&\sum_{\beta_h>0} \frac{\beta_h}{\|z_{p(h),h}\|} x_{p(h),h}+\sum_{\beta_h<0} \frac{-\beta_h}{\|z_{p(h+d-1),h+d-1}\|}x_{p(h+d-1),h+d-1}+\mu v\\
\end{eqnarray*}
where we collected the $\lambda_{p(h),h}v$ terms in the single term $\mu v$. Then $a=(b-\mu v)+(\lambda +\mu) v$. The previous formula shows that $b-\mu v$ is a non-negative combination of the elements in $T$. The term $(\lambda +\mu) v$ or $(-\lambda -\mu) (-v)$ is taken as a non-negative multiple of $v\in X_{2d-1}$ or of $-v\in X_{2d}$ depending on the sign of $\lambda+\mu$. \qed

\smallskip
The uniqueness part of Theorem~\ref{th:stein}, $\pos T^*=\R^d$ and our assumption imply that, for every permutation $p$, if $T$ contains the vector $u$ then it contains $-u$. Assume for instance that $x_{1,1}=u$ and $p(1)=1$ and $p(d)=j$ (of course $j\ne 1$). Then $x_{j,d}=x_{p(d),d}=-u$ and no other $x_{p(h),h}\in T$ can be equal to $-u$. This implies $x_{j,d}=-u$ for all $j\in [2d-2],j\ne 1$, and a simple argument shows that $x_{1,d}=-u$ as well. Thus every entry in column $d$ of $M(\X^{\circ})$ is $-u=-x_{1,1}$. 

\smallskip
This proof works for $x_{1,i}$ for every $i<d$ and shows that every entry in column $i+d-1$ of $M(\X^{\circ})$ is equal to $-x_{1,i}$. Symmetrically, every entry in column $i< d$ is equal to $-x_{1,i+d-1}$. Setting now $x_{1,i}=e_i$ for every $i<d$ finishes the proof.\qed

\bigskip
\section{Proof of Lemma~\ref{l:subPC}}\label{sec:subPC}

The proof is very similar to the previous one so we only give a sketch. The matrix representation of $\W$ is a $(2d-2)\times d$ matrix $(w_{i,j})$ where $w_{i,j}=f_j$ when $i<d$ and 
$w_{i,j}=-f_j$ when $i\ge d$. This time a transversal $T$ of $\X^{\circ}$ is given by $k \in [d]$ and two one-to-one maps $p^+$ and $p^-$. 
We extend it again to the transversal $T^*$ of $\X$ by adding $v\in X_{2d-1}$ and $-v\in X_{2d}$.

\begin{claim} $\pos T^*=\R^d$.
\end{claim}

The {\bf proof} uses equation (\ref{eq:bPos}) and is similar to the previous one and is omitted. 

\smallskip
We conclude again that if $T$ contains a vector $u$ then it contains $-u$ as well. This time the first $d-1$ entries in column $j$ of $M(\X^{\circ})$ are equal to $x_{1,j}$ and the rest to  $-x_{1,j}$. Defining $f_j=x_{1,j}$ for $j\in [d]$ completes the proof.\qed

\bigskip
{\bf Acknowledgements.} The first author (IB) was partially supported by NKFIH grant No. 133819 and also by the HUN-REN Research Network. 
The second author (YQ) was supported by the China Scholarship Council (grant No. 202508810002) and the 2026 Graduate Student Innovation Funding Project under the College of Hebei Normal University (grant No. ycxzzbs202604). This work was additionally supported by the following grants: the NSF of China (122711392); the Foreign Experts Program of the People’s Republic of China; the Program for Foreign Experts of Hebei Province.

\vskip0.6cm

\newpage
\noindent
Imre B\'ar\'any \\

\noindent
Alfr\'ed R\'enyi Institute of Mathematics, HUN-REN\\
13 Re\'altanoda Street, Budapest 1053 Hungary,\\
{\tt barany.imre@renyi.hu} \\

\noindent
Department of Mathematics, University College London\\
Gower Street, London, WC1E 6BT, UK, and \\

\noindent
School of Mathematical Sciences,
Hebei Normal University,
050024 Shijiazhuang, P.R. China,\\

\noindent
Yun Qi\\

\noindent
School of Mathematical Sciences,
Hebei Normal University,
050024 Shijiazhuang, P.R. China.\\
{\tt yunqi1632024@163.com}
\end{document}